\documentclass[10pt,twoside]{article}
\usepackage{graphicx}
\usepackage{amsmath}
\usepackage{Latex-document}

\title{\bf  Motivic Integration and the Grothendieck Group of Pseudo-Finite Fields \vskip 6mm}

\author{J. Denef\,\vspace*{-0.5cm}\thanks{Department of Mathematics,
University of Leuven, Celestijnenlaan 200 B, 3001 Leuven,
Belgium. E-mail: Jan.Denef@wis.kuleuven.ac.be}\; and F. Loeser\vspace*{-0.5cm}\thanks{D\'epartement
de Math\'ematiques et Applications,
\'Ecole Normale Sup\'erieure, 45 rue d'Ulm, 75230 Paris Cedex 05, France (UMR 8553 du
CNRS). E-mail: Francois.Loeser@ens.fr}  }

\date{}

\begin{document}
\maketitle

\thispagestyle{first} \setcounter{page}{1}

\begin{abstract}\vskip 3mm

Motivic integration is a powerful
technique to prove that certain quantities associated to algebraic
varieties are birational invariants or are independent of a chosen resolution of singularities.
We survey our recent work on an extension of the
theory of motivic integration,  called arithmetic motivic integration. We developed this
theory to understand how $p$-adic integrals of a very general type depend on $p$.
Quantifier elimination plays a key role.
\vskip 4.5mm


\end{abstract}

\vskip 12mm

\section{Introduction} \label{section 1}
\setzero \vskip-5mm \hspace{5mm }

Motivic integration was first introduced  by Kontsevich \cite{Kontsevich} and further developed by
Batyrev \cite{Batyrev: Stringy Hodge}\cite{Batyrev: McKay}, and Denef-Loeser
\cite{Denef Loeser: Inventiones}\cite{Denef Loeser: Thom Sebastiani}\cite{Denef Loeser: McKay}.
It is a powerful
technique to prove that certain quantities associated to algebraic
varieties are birational invariants or are independent of a chosen resolution of singularities.
For example, Kontsevich used it to prove that the Hodge numbers of birationally equivalent
projective Calabi-Yau manifolds are equal. Batyrev \cite{Batyrev: Stringy Hodge} obtained  his string-theoretic  Hodge numbers
for canonical Gorenstein singularities by motivic integration. These are the right quantities to
establish several mirror-symmetry identities for Calabi-Yau varieties.
For more applications and references we refer to the survey
papers \cite{Denef Loeser: Barcelona} and \cite{Looijenga}. Since
than, several other applications to singularity theory were discovered,
see e.g. Musta\c t\u a \cite{Mustata: Jet Schemes}.

In the present paper, we survey our recent work \cite{Denef Loeser: Definable Sets} on an extension of the
theory of motivic integration,  called arithmetic motivic integration. We developed this
theory to understand how $p$-adic integrals of a very general type depend on $p$. This is
used in recent work of  Hales \cite{Hales} on orbital integrals related to the Langlands program. Arithmetic
motivic integration is tightly linked to the theory of quantifier elimination, a subject
belonging to mathematical logic. The roots of this subject go back to Tarski's theorem on
projections of semi-algebraic sets and to the  work of Ax-Kochen-Ersov and Macintyre on
quantifier elimination for Henselian valued fields (cf. section 4). We will illustrate arithmetic
motivic integration starting with the following concrete application. Let $X$ be an algebraic
variety given
by equations with integer coefficients. Denote by $ N_{p,n}$ the cardinality of the image of the
projection $ X({\bf Z}_p ) \to X({\bf Z}/p^{n+1} )$, where ${\bf Z}_p $ denotes the  $p$-adic
integers. A conjecture of Serre and Oesterl\'e states that
$P_p (T): = \sum\limits_n {N_{p,n} } T^n$ is rational. This was proved in 1983
by Denef \cite{Denef: Inventiones 83} using
quantifier elimination, expressing $P_p (T)$ as a $p$-adic integral over a domain defined by a
formula involving quantifiers. This gave no information yet on how $P_p (T)$ depends
on $p$. But recently, using arithmetic motivic integration, we  proved:

{\bf Theorem 1.1 } \it
There exists a canonically defined rational power series $P(T)$ over the ring
${\rm K}_0^{mot} ({\rm Var}_{\bf Q}  ) \otimes {\bf Q}$, such that, for $p\gg$ 0,  $P_p(T)$ is
obtained from $P(T)$ by applying to each coefficient of $P(T)$ the operator $N_p$.
\rm

Here ${\rm K}_0 ({\rm Var}_{\bf Q}  )$
denotes the Grothendieck ring of algebraic varieties over $\bf Q$, and $
{\rm K}_0^{mot} ({\rm Var}_{\bf Q}  )$ is the quotient of this ring obtained by identifying two varieties
if they have the same class in the Grothendieck group of Chow
motives (this is explained in the next section). Moreover the operator
$N_p$ is induced by  associating to a variety over $\bf Q$ its number of rational points over
the field with $p$ elements, for $p\gg 0$.

As explained in section 8 below, this theorem is a special case of a
much more general theorem on $p$-adic
integrals. There we will also see how to canonically associate
a \textquotedblleft virtual motive\textquotedblright\: to quite general $p$-adic integrals.
A first step in the proof of the above theorem is the construction of a
canonical morphism from the Grothendieck ring ${\rm K}_0 ({\rm PFF}_{\bf Q} )$ of the theory of
pseudo-finite fields of characteristic zero, to ${\rm K}_0^{mot} ({\rm Var}_{\bf Q}  ) \otimes {\bf Q}$.
Pseudo-finite fields play a key role in the work of Ax \cite{Ax} that leads to  quantifier
elimination for finite fields \cite{Kiefe}\cite{Fried Sacerdote}\cite{Chatzidakis}. The existence of this map is interesting in
itself, because any generalized Euler characteristic, such as the topological Euler
characteristic or the Hodge-Deligne polynomial, can be evaluated on any element of
${\rm K}_0^{mot} ({\rm Var}_{\bf Q}  ) \otimes {\bf Q}$, and hence also on any logical formula in the
language of fields (possibly involving quantifiers). All this will be explained in
section 2. In section 3 we state Theorem 3.1, which is a stronger version of Theorem 1.1
that determines $P(T)$. A proof of Theorem 3.1 is outlined in section 7, after giving a
survey on arithmetic motivic integration in section 6.

\section{The Grothendieck group of pseudo-finite fields}
\label{section 2} \setzero\vskip-5mm \hspace{5mm }

Let $k$ be a field of characteristic zero. We denote by ${\rm K}_0 ({\rm Var}_{k}  )$
the Grothendieck ring of algebraic varieties over $k$. This is the group generated by
symbols $[V]$ with $ V $ an algebraic variety over $k$, subject to the relations
$[V_1]=[V_2]$ if $V_1$ is isomorphic to $V_2$, and $[V\setminus W]= [V]-[W]$ if $W$
is a Zariski closed subvariety of $W$.
The ring multiplication on ${\rm K}_0 ({\rm Var}_{k}  )$ is induced by the cartesian
product of varieties. Let {\bf L} be the class of the affine line over $k$
in ${\rm K}_0 ({\rm Var}_{k}  )$. When $V$ is an algebraic variety over {\bf Q},
and $p$ a prime number, we
denote by $N_p(V)$ the number of rational points over the field ${\bf F}_p$ with $p$ elements on a
model $\tilde{V}$ of $V$ over {\bf Z}. This depends on the choice of a model $\tilde{V}$,
but two different models will yield the same value of $N_p(V)$, when $p$ is large enough.
This will not cause any abuse later on. For us, an algebraic variety over $k$ does not need to be
irreducible; we mean by it a reduced separated scheme of finite over $k$.

To any projective nonsingular variety over $k$ one
associates its Chow motive over $k$ (see \cite{Scholl}). This is
a purely algebro-geometric construction, which is
made in such a way that any two projective nonsingular varieties, $ V_1$ and $V_2$, with isomorphic
associated Chow motives, have the same cohomology for each of the known cohomology
theories (with coefficients in a field of characteristic zero). In particular, when $k$ is {\bf Q}, $N_p (V_1)=N_p (V_2)$, for $p \gg 0$.
For example two elliptic curves define the same Chow motive iff there is a
surjective morphism from one to the other.
We denote by $ {\rm K}_0^{mot} ({\rm Var}_{k}  )$ the quotient of the ring ${\rm K}_0 ({\rm Var}_{k}  )$
obtained by identifying any two nonsingular projective varieties over $k$ with equal
associated Chow motives.  From work of  Gillet
and Soul\'e \cite{Gillet Soule}, and Guill\' en and Navarro Aznar\cite{Guillen Aznar}, it
directly follows that there is a unique ring monomorphism from $ {\rm K}_0^{mot} ({\rm Var}_{k}  )$
to the Grothendieck
ring of the category of Chow motives over $k$, that
maps the class of a projective nonsingular variety to the class of its associated Chow motive.
What is important for the applications, is that any generalized
Euler characteristic, which can be defined in terms of cohomology (with coefficients in a field
of characteristic zero),
factors through $ {\rm K}_0^{mot} ({\rm Var}_{k}  )$. With a generalized
Euler characteristic we mean any ring morphism from ${\rm K}_0 ({\rm Var}_{k}  )$, for
example the topological  Euler characteristic and the Hodge-Deligne polynomial when $k$ = {\bf
C}.
For $[V]$ in $ {\rm K}_0^{mot} ({\rm Var}_{k}  )$, with $k$ = {\bf
Q}, we put $N_p ([V])= N_p(V)$; here again
this depends on choices, but two different choices yield the same value for $N_p ([V])$, when
$p$ is large enough.

With a ring formula $\varphi$ over $k$ we
mean a logical formula build from polynomial equations over $k$, by taking Boolean
combinations and using existential and universal quantifiers.
For example, $(\exists x)(x^2+x+y=0 \; {\rm and}\; 4y\neq 1)$
is a ring formula over {\bf Q}.
The mean purpose of the
present section is to associate in a canonical way to each such formula $\varphi$ an
element $\chi _c ([\varphi])$
of ${\rm K}_0^{mot} ({\rm Var}_{k}  ) \otimes {\bf Q} $. One of the required
properties of this association is
the following, when $k$ = {\bf Q}: If the formulas $\varphi_1$ and $\varphi_2$ are
equivalent when interpreted in ${\bf F}_p$, for all large enough primes $p$, then
$\chi _c ([\varphi_1])$ = $\chi _c ([\varphi_2])$. The natural generalization of this
requirement, to arbitrary fields $k$ of characteristic zero, is the following:
If the formulas $\varphi_1$ and $\varphi_2$ are
equivalent when interpreted in $K$, for all pseudo-finite fields $K$ containing $k$, then
$\chi _c ([\varphi_1])$ = $\chi _c ([\varphi_2])$. We recall that a
pseudo-finite field is an infinite perfect field that has exactly one field extension of any given finite
degree, and over which each absolutely irreducible variety has a rational point.
For example, infinite ultraproducts of finite fields are pseudo-finite.
J. Ax \cite{Ax}
proved that two ring formulas over {\bf Q} are equivalent
when interpreted in ${\bf F}_p$, for all large enough primes $p$, if and only if they are
equivalent when interpreted in $K$, for all pseudo-finite fields $K$ containing ${\bf
Q}$. This shows that the two above mentioned requirements are equivalent when $k$ = {\bf
Q}. In fact, we will require much more, namely that the association
$\varphi$ $\longmapsto$ $\chi _c ([\varphi])$ factors through the Grothendieck ring
${\rm K}_0 ({\rm PFF}_k )$ of the theory of pseudo-finite fields containing $k$. This
ring is the group generated by symbols $[\varphi ]$, where $\varphi$ is any ring formula
over $k$, subject to the relations $[\varphi_1 \,{\rm or}\;
\varphi_2]=[\varphi_1]+[\varphi_2]- [\varphi_1 \, {\rm and}\;
\varphi_2]$, whenever $\varphi_1$ and $\varphi_2$ have the same free variables, and the
relations
$[\varphi_1] = [\varphi_2]$, whenever there exists a
ring formula $\psi$ over $k$ that, when interpreted in any pseudo-finite field $K$ containing $k$,
yields the graph of a bijection between the tuples of elements of $K$ satisfying $\varphi_1$
and those satisfying $\varphi_2$. The ring multiplication on ${\rm K}_0 ({\rm PFF}_k )$
is induced by the conjunction of formulas in disjoint sets of variables. We can now state
the following variant of a theorem of Denef and Loeser \cite{Denef Loeser: Definable Sets}.

{\bf Theorem 2.1.} \it There exists a unique ring morphism $$ \chi _c : {\rm K}_0
({\rm PFF}_k ) \longrightarrow {\rm K}_0^{mot} ({\rm Var}_{k}  ) \otimes {\bf Q} $$
satisfying  the following two properties:

\noindent (i) For any formula $\varphi$ which is a conjunction of polynomial equations
over $k$, the element $\chi _c ([\varphi])$ equals the
class in ${\rm K}_0^{mot} ({\rm Var}_{k}  ) \otimes {\bf Q}$ of the variety defined
by $\varphi$.

\noindent (ii) Let $X$ be a normal affine irreducible variety over $k$, $Y$ an unramified
Galois cover \footnote{Meaning that $Y$ is an integral \'etale scheme over
$X$ with $Y/G \cong X$, where $G$ is the group of all endomorphisms of $Y$ over $X$.}  of
$X$, and $C$ a cyclic subgroup of the Galois group G of $Y$ over $X$. For such data we
denote by $\varphi _{Y,X,C}$ a ring formula, whose interpretation in any field $K$ containing
$k$, is the set of $K$-rational points on $X$ that lift to a geometric
point on $Y$ with decomposition group $C$ (i.e. the set of points on $X$ that lift to a $K$-rational
point of $Y/C$, but not to any $K$-rational point of $Y/C'$ with $C'$ a proper subgroup of
$C$). Then
\[
\chi _c ([\varphi _{Y,X,C} ]) = \frac{{\left| C \right|}}{{\left| {{\rm N}_{\rm G} (C)} \right|}}\chi _c ([\varphi _{Y,Y/C,C}
]),\]
where ${\rm N}_{\rm G} (C)$ is the normalizer of $C$ in $G$.

Moreover, when  $k$ = {\bf Q}, we have for all large enough primes $p$ that
$ N_p (\chi _c ([\varphi]))$ equals the number of tuples in ${\bf F}_p$ that satisfy the
interpretation of $\varphi$ in ${\bf F}_p$. \rm

The proof of the uniqueness goes as follows: From quantifier elimination for
pseudo-finite fields (in terms of Galois stratifications,
cf. the work of Fried and Sacerdote \cite{Fried Sacerdote}\cite[\S 26]{Fried Jarden}), it follows that every ring
formula over $k$ is equivalent (in all pseudo-finite fields containing $k$) to a Boolean
combination of formulas of the form $\varphi _{Y,X,C}$. Thus by (ii) we only have to
determine $\chi_c([\varphi _{Y,Y/C,C}])$, with $C$ a cyclic group. But this follows
directly from the following recursion formula: \[
\left| C \right|[Y/C] = \sum\limits_{A\,{\rm  subgroup\,  of }\, C} {\left| A \right|} \chi _c ([\varphi _{Y,Y/A,A}
]).\]
This recursion formula is a direct consequence of (i), (ii), and the fact that
the formulas $\varphi _{Y,Y/C,A}$ yield a partition of $Y/C$.
 The proof of the existence
of the morphism $\chi _c$ is based on the following. In \cite{Rollin Aznar}, del Ba{\~ n}o
Rollin and
Navarro Aznar associate to any representation over {\bf Q} of a finite group $G$ acting
freely on an affine variety $Y$ over $k$, an element in the Grothendieck group of Chow
motives over $k$. By linearity, we can hence associate to any {\bf Q}-central function
$\alpha$ on $G$ (i.e. a {\bf Q}-linear combination of characters of
representations of $G$ over {\bf Q}), an element $\chi_c(Y,\alpha )$ of that Grothendieck
group tensored with {\bf Q}. Using Emil Artin's Theorem, that any {\bf Q}-central function
$\alpha$ on $G$ is a {\bf Q}-linear combination of characters induced by trivial
representations of cyclic subgroups, one shows that
$\chi_c(Y,\alpha )\in {\rm K}_0^{mot} ({\rm Var}_{k}  ) \otimes {\bf Q}$. For $X:=Y/G$
and $C$ any cyclic subgroup of $G$,
we define $\chi_c([\varphi _{Y,X,C}]) := \chi_c(Y,\theta)$, where $\theta$ sends $g\in G$ to 1
if the subgroup generated by $g$ is conjugate to $C$, and else to 0. Note that
$\theta$ equals
$\left| C \right|/\left| {N_G (C)} \right|$ times the function on G induced by
the characteristic function on $C$ of the set of generators of $C$. This implies our
requirement (ii), because of Proposition 3.1.2.(2) of \cite{Denef Loeser: Definable Sets}.
The map
$(Y,\alpha)\mapsto \chi_c (Y,\alpha)$ satisfies the nice compatibility relations
stated in Proposition 3.1.2 of loc. cit. This compatibility (together with the above
mentioned quantifier elimination) is used, exactly as in loc. cit., to prove that
the above definition of $\chi_c([\varphi _{Y,X,C}])$ extends by additivity to a
well-defined map
$ \chi _c : {\rm K}_0
({\rm PFF}_k ) \longrightarrow {\rm K}_0^{mot} ({\rm Var}_{k}  ) \otimes {\bf Q} $.
In loc. cit., Chow motives with coefficients in the algebraic closure of {\bf Q} are used,
but we can work as well with coefficients in {\bf Q}, since here we only have to consider
representations of $G$ over {\bf Q}.

\section{Arc spaces and the motivic Poincar\'e series}
\label{section 3} \setzero\vskip-5mm \hspace{5mm }

Let $X$ be an algebraic variety defined over a field $k$ of characteristic zero. For any
natural number $n$, the $n$-th jet space $\mathcal{L}_n (X)$ of $X$ is the unique algebraic
variety over $k$ whose $K$-rational points correspond in a bijective and functorial way to the
rational points on $X$ over $K[t]/t^{n+1}$, for any field $K$ containing $k$. The arc
space $\mathcal{L} (X)$ of $X$ is the reduced $k$-scheme obtained by taking the
projective limit of the varieties $\mathcal{L}_n (X)$ in the category of $k$-schemes.

We will now give the definition of the motivic Poincar\'e series $P(T)$ of $X$.
This series is called the arithmetic Poincar\'e  series in \cite{Denef Loeser: Definable
Sets}, and is very different from the geometric Poincar\'e  series studied
in \cite{Denef Loeser: Inventiones}.
For notational convenience we only give the definition here when $X$ is a subvariety of
some affine space $\mathbf{A}_k^m$. For the general case we refer to section 5 below or to
our paper \cite{Denef Loeser: Definable Sets}.
By Greenberg's Theorem \cite{Greenberg}, for each $n$ there exists a ring formula $\varphi_n$ over $k$  such
that, for all fields $K$ containing $k$, the $K$-rational points of
$\mathcal{L}_n (X)$, that can be lifted to a $K$-rational point of $\mathcal{L} (X)$,
correspond  to the tuples satisfying the interpretation of $\varphi_n$ in $K$. (The
correspondence is induced by mapping a polynomial over $K$ to the tuple consisting of its
coefficients.) Clearly, when two formulas satisfy this requirement, then
they are equivalent when interpreted in any field containing $k$, and hence define the same
class in ${\rm K}_0 ({\rm PFF}_k )$. Now we are ready to give the definition of $P(T)$:
\[
P(T): = \sum\limits_n {\chi _c ([\varphi _n ])T^n }
\]

{\bf Theorem 3.1.} \it The motivic Poincar\'e series $P(T)$ is a rational power series
over the ring ${\rm K}_0^{mot} ({\rm Var}_{k}  )[{\bf L}^{-1}] \otimes {\bf Q}$, with denominator
a product of factors of the form $1-{\bf L}^{a}T^b$, with $a,b \in{\bf Z}$,
$b>0$. Moreover if $k={\bf Q}$, the Serre Poincar\'e series
$P_p(T)$, for $p\gg$0, is obtained
from $P(T)$ by applying the operator $N_p$ to each coefficient of the numerator and
denominator of $P(T)$.
\rm

In particular we see that the degrees of the numerator and the denominator of $P_p(T)$
remain bounded for $p$ going to infinity. This fact was first proved by Macintyre
\cite{Macintyre: bounded} and Pas \cite{Pas}.

\section{Quantifier elimination for valuation rings}
\label{section 4} \setzero\vskip-5mm \hspace{5mm }

Let $R$ be a ring and assume it is an integral domain. We will define the notion of a DVR-formula
over $R$. Such a formula can be
interpreted in any discrete valuation ring $A \supset R$ with a distinguished uniformizer $\pi$.
It can contain variables that run over the discrete valuation ring, variables that run
over the value group {\bf Z}, and variables that run over the residue field. A DVR-formula
over $R$ is build from quantifiers with respect to variables that run over the discrete
valuation ring, or over the value group, or over the residue field, Boolean combinations,
and expressions of the following form: $g_1(x)=0$,
${\rm ord}(g_1(x))\leq {\rm ord}(g_2(x))+ L(a)$,
${\rm ord}(g_1(x))\equiv L(a) \; {\rm mod} \; d $,
where $g_1(x)$ and $g_2(x)$ are polynomials
over $R$ in several variables $x$ running over the discrete valuation ring, where $L(a)$
is a polynomial of degree $\leq 1$ over {\bf Z} in several variables $a$ running
over the value group, and $d$ is any positive integer (not a variable). Moreover we also allow
expressions of the
form $\varphi (\overline{{\rm ac}}(h_1(x)),...,\overline{{\rm ac}}(h_r(x)))$, where
$\varphi$
is a ring formula over $R$, to be interpreted in the residue field,
$h_1(x),...,h_r(x)$ are polynomials over $R$ in several variables $x$ running over the
discrete valuation ring, and $\overline{{\rm ac}}(v)$, for any element $v$ of the discrete
valuation ring, is the residue of the angular
component ${\rm ac}(v):= v \pi ^{-{\rm ord  }v}$. For the discrete valuation rings ${\bf Z}_p$
and $K[[t]]$, we take as distinguished uniformizer $\pi $ the elements $p$ and $t$.

{\bf Theorem 4.1 \ (Quantifier Elimination of Pas \cite{Pas})} \it Suppose that $R$ has
characteristic zero. For any DVR-formula $\theta$ over $R$ there exists a DVR-formula
$\psi$ over $R$, which contains no quantifiers running over the valuation ring and no
quantifiers running over the value group, such that \\
(1) $\theta \longleftrightarrow \psi$ holds in $K[[t]]$, for all fields $K$ containing
$R$. \\
(2) $\theta \longleftrightarrow \psi$ holds in ${\bf Z}_p$, for all primes $p\gg 0$, when
$R= {\bf Z}$.
\rm

The Theorem of Pas is one of several quantifier elimination results for Henselian valuation
rings, and goes back to the work of Ax-Kochen-Ersov and Cohen on the model theory
of valued fields,
which was further developed by Macintyre, Delon \cite{Delon}, and others, see e.g. Macintyre's
survey \cite{Macintyre: survey}.

Combining the
Theorem of Pas with the work of Ax mentioned in section 2, one obtains

{\bf Theorem 4.2 \ (Ax-Kochen-Ersov Principle, version of Pas)} \it Let $\sigma$ be a
DVR-formula over {\bf Z} with no free variables. Then the following are equivalent: \\
(i)  The interpretation of $\sigma$ in ${\bf Z}_p$ is true for all primes $p\gg 0$. \\
(ii) The interpretation of $\sigma$ in $K[[t]]$ is true for all pseudo-finite fields $K$
of characteristic zero.

\rm

\section{Definable subassignements and truncations}
\label{section 5} \setzero\vskip-5mm \hspace{5mm }

Let $h : \mathcal{C} \rightarrow {\rm Sets}$ be a functor from a
category $\mathcal{C}$ to the category of sets. We shall call
the data for each object $C$ of $\mathcal{C}$ of a subset $h' (C)$ of
$h (C)$  a {\em subassignement} of $h$. The point in this definition
is that
$h'$ is not assumed to
be a subfunctor of $h$.
For $h'$ and $h''$ two subassignements
of $h$, we shall denote by $h' \cap h''$ and $h' \cup h''$,
the subassignements
$ C \mapsto h' (C) \cap h'' (C)$ and $ C \mapsto h' (C) \cup h'' (C)$,
respectively.

Let $k$ be a field of characteristic zero. We denote by ${\rm Field}_{k}$ the category
of fields which contain $k$.
For $X$ a variety over $k$, we consider the functor $h_{X} : K \mapsto X
(K)$ from ${\rm Field}_{k}$ to the category of sets. Here $X(K)$ denotes the set of
$K$-rational points on $X$. When $X$ is a subvariety of some affine space, then a
subassignement $h$ of $h_{X}$ is called {\em definable} if there exists a ring  formula $\varphi$ over $k$
such that, for any field $K$ containing $k$, the set of tuples that satisfy the
interpretation of $\varphi$ in $K$, equals $h(K)$. Moreover we define the {\em class $[h]$ of $h$}
in ${\rm K}_0 ({\rm PFF}_k )$ as $[\varphi]$.
More generally, for any algebraic variety $X$ over $k$, a subassignement $h$ of $h_{X}$
is called {\em definable} if  there exists a finite cover
$(X_{i})_{i \in I}$ of $X$ by affine open subvarieties and
definable subassignements $h_{i}$ of $h_{X_{i}}$, for $i \in I$,
such that
$h  = \cup_{i \in I} h_{i}$. The {\em class $[h]$ of $h$}
in ${\rm K}_0 ({\rm PFF}_k )$ is defined by linearity, reducing to the affine case.

For any algebraic variety $X$ over $k$ we denote by $h_{\mathcal{L} (X)}$ the functor
$ K \mapsto X(K[[t]])$ from ${\rm Field}_{k}$ to the category of sets. Here $X(K[[t]])$
denotes the set of $K[[t]]$-rational points on $X$. When $X$ is a subvariety of some affine space,
then a subassignement $h$ of $h_{\mathcal{L} (X)}$ is called {\em definable} if there exists
a DVR-formula $\varphi$ over $k$
such that, for any field $K$ containing $k$, the set of tuples that satisfy the
interpretation of $\varphi$ in $K[[t]]$, equals $h(K)$.
More generally, for any algebraic variety $X$ over $k$, a
subassignement $h$ of $h_{\mathcal{L} (X)}$
is called {\em definable} if  there exists a finite cover
$(X_{i})_{i \in I}$ of $X$ by affine open subvarieties and
definable subassignements $h_{i}$ of $h_{\mathcal{L} (X_{i})}$, for $i \in I$,
such that
$h  = \cup_{i \in I} h_{i}$.
A family of definable subassignements $h_n$, $n\in {\bf Z}$, of $h_{\mathcal{L} (X)}$ is
called a {\em definable family of definable subassignements} if on each affine open of a
suitable finite affine covering of $X$, the family $h_n$ is given by a DVR-formula
containing $n$ as a free variable running over the value group.

Let $X$ be a variety over $k$. Let $h$ be a definable subassignement of $h_{\mathcal{L}
(X)}$, and $n$ a natural number. The {\em truncation of $h$ at level $n$, denoted by
$\pi_n (h)$}, is the subassignement of $h_{\mathcal{L}_n (X)}$ that associates to any field
$K$ containing $k$ the image of $h(K)$ under the natural projection map from
$X(K[[t]])$ to $\mathcal{L}_n (X)(K)$. Using the Quantifier Elimination Theorem of Pas, we proved
that $\pi_n (h)$ is a definable subassignement of $h_{\mathcal{L}_n (X)}$, so that we can
consider its class [$\pi_n (h)$] in ${\rm K}_0 ({\rm PFF}_k )$. Using the notion of
truncations, we can now give an alternative (but equivalent) definition of the motivic
Poincar\'e series $P(T)$, which works for any algebraic variety $X$ over $k$, namely
$P(T): = \sum\limits_n {\chi _c ([\pi_n ( h_{\mathcal{L} (X)} )])T^n }$.

A  definable subassignement $h$ of $h_{\mathcal{L}(X)}$ is called
\emph{weakly stable at level $n$ } if
for any field $K$ containing $k$ the set $h(K)$ is a union of fibers of the natural
projection map from $X(K[[t]])$ to $\mathcal{L}_n (X)(K)$. If $X$ is nonsingular, with
all its irreducible components of dimension $d$, and $h$ is a
definable subassignement of $h_{\mathcal{L}(X)}$, which is weakly stable at level $n$, then
it is easy to verify that \[
[\pi _n (h)]{\bf L}^{ - nd}  = [\pi _m (h)]{\bf L}^{ - md}
\]
for all $m\geq n$. Indeed this follows from the fact that the natural map from
$\mathcal{L}_m (X)$ to $\mathcal{L}_n (X)$ is a locally trivial fibration for the Zariski
topology with fiber $\mathbf{A}_k^{(m-n)d}$, when $X$ is nonsingular.

\section{Arithmetic motivic integration}
\label{section 6} \setzero\vskip-5mm \hspace{5mm }

Here we will outline an extension of the
theory of motivic integration, called arithmetic motivic integration. If the base field
$k$ is algebraically closed, then it coincides with the usual motivic integration.

 We denote by
$\widehat{\rm K}_0^{mot} ({\rm Var}_k )[{\bf L}^{ - {\bf 1}} ]$
the completion of
${\rm K}_0^{mot} ({\rm Var}_k )[{\bf L}^{ - {\bf 1}} ]$ with respect to the filtration of
${\rm K}_0^{mot} ({\rm Var}_k )[{\bf L}^{ - {\bf 1}} ]$ whose $m$-th member is the
subgroup generated by the elements $[V]{\bf L}^{-i}$ with $i-{\rm dim}V \geq m$.
Thus a sequence $[V_i]{\bf L}^{-i}$ converges to zero in
$\widehat{\rm K}_0^{mot} ({\rm Var}_k )[{\bf L}^{ - {\bf 1}} ]$, for $i \mapsto +\infty$,
if $i-{\rm dim}V_i$  $\mapsto +\infty$.

{\bf Definition-Theorem 6.1.} \it Let $X$ be an algebraic variety of dimension $d$
over a field $k$ of characteristic zero, and let $h$ be a definable subassignement
of $h_{\mathcal{L}(X)}$. Then the limit \[ \nu (h):=
\mathop {{\rm lim}}\limits_{n \to \infty }\chi_c([\pi _n (h)]){\bf L}^{ - (n+1)d}
\] exists in $\widehat{\rm K}_0^{mot} ({\rm Var}_k )[{\bf L}^{ -  1} ]\otimes {\bf Q }$ and is
called the arithmetic motivic volume of $h$.

\rm
We refer to \cite[\S 6]{Denef Loeser: Definable Sets}  for the proof of the above theorem. If $X$ is
nonsingular and $h$ is weakly stable at some level, then the theorem follows directly
from what we said at the end of the previous section. When $X$ is nonsingular affine, but $h$
general, the theorem is proved by approximating $h$ by definable subassignements $h_i$
of $h_{\mathcal{L}(X)}$, $i\in {\bf N}$, which are weakly stable at level $n(i)$.
For $h_i$ we take the subassignement obtained from $h$ by adding, in the
DVR-formula $\varphi$ defining $h$, the condition
$ {\rm ord }  g(x) \leq i$, for each polynomial $g(x)$ over the valuation ring,
that appears in $\varphi$.
(Here we assume that $\varphi$ contains no
quantifiers over the valuation ring.) It
remains to show then that $\chi_c([\pi_n({\rm ord } g(x) > i)]){\bf L}^{ - (n+1)d}$ goes to zero
when both $i$ and $n\gg i$ go to infinity, but this is easy.

{\bf Theorem 6.2.} \it Let $X$ be an algebraic variety of dimension $d$
over a field $k$ of characteristic zero, and let $h$, $h_1$ and $h_2$ be definable
subassignements of $h_{\mathcal{L}(X)}$.

\noindent (1) If $h_1(K)=h_2(K)$ for any pseudo-finite field $K \supset k$, then
$\nu (h_1)=\nu (h_2)$.

\noindent (2) $ \nu (h_1  \cup h_2)= \nu (h_1)+ \nu (h_2)- \nu (h_1 \cap h_2)$

\noindent (3) If $S$ is a  subvariety of $X$ of dimension $<d$, and if
$h \subset h_{\mathcal{L}(S)}$, then $\nu (h)=0$.

\noindent (4) Let $h_n$, $n\in {\bf N }$, be a definable family of definable
subassignements of $h_{\mathcal{L}(X)}$. If $h_n \cap h_m =\emptyset$, for all $n \neq
m$, then $\sum\limits_n {\nu (h_n )}$ is convergent and equals $\nu( \bigcup \limits_n
{h_n)}$.

\noindent (5) \rm  Change of variables formula. \it Let $p:Y \to X$ be a proper
birational morphism of nonsingular irreducible
varieties over $k$. Assume for any field $K$ containing $k$ that the jacobian determinant
of $p$ at any point of $p^{-1}(h(K))$ in $ Y (K[[t]])$ has t-order equal to $e$. Then $\nu_X
(h)= {\bf L}^{-e} \nu_Y (p^{-1}(h))$. Here $\nu_X$, $\nu_Y$ denote the arithmetic motivic
volumes relative to $X$, $Y$, and $p^{-1}(h)$ is the subassignement of $h_{\mathcal{L}(Y)}$
given by $K \mapsto p^{-1}(h(K)) \cap Y (K[[t]])$.
\rm

Assertion (1) is a direct consequence of the definitions. Assertions (2) and (4) are
proved by approximating the subassignements by weakly stable ones. Moreover for (4) we
also need the fact that $h_n= \emptyset$ for all but a finite number of $n$'s, when all
the $h_n$, and their union, are weakly stable (at some level depending on $n$). Assertion (5) follows from
the fact that for $n \gg e$ the map $\mathcal{L}_n (Y)   \rightarrow \mathcal{L}_n (X)$
induced by $p$ is a piecewise trivial fibration with fiber ${\bf A }_k^e$ over the image
in $\mathcal{L}_n (X)$ of the points of $\mathcal{L} (Y)$ where the jacobian determinant
of $p$ has $t$-order $e$. See \cite{Denef Loeser: Definable Sets} for the details.

\section{About the proof of Theorem 3.1}
\label{section 7} \setzero\vskip-5mm \hspace{5mm }

We give a brief sketch of the proof of Theorem 3.1, in the special case that $X$ is a
hypersurface in ${\bf A}_k ^d$ with equation $f(x)=0$. Actually, here we will only
explain why the image $\widehat{P}(T)$ of $P(T)$ in the ring of power series over
$\widehat{\rm K}_0^{mot} ({\rm Var}_k )[{\bf L}^{ -  1} ]\otimes {\bf Q }$ is rational.
The rationality of $P(T)$ requires additional work.
Let $\varphi (x,n)$ be the
DVR-formula $(\exists y)(f(y)=0 \;{\rm and}\;{\rm ord} (x-y)\geq n)$, with $d$ free
variables $x$ running over the discrete valuation ring, and one free variable $n$ running
over the value group. That formula determines a definable family of definable
subassignements $h_{\varphi (-,n)}$ of $h_{\mathcal{L}({\bf A}_k ^d)}$. Since
$h_{\varphi (-,n)}$ is weakly stable at level $n$, unwinding our definitions yields that
the arithmetic motivic volume on
$h_{\mathcal{L}({\bf A}_k ^d)}$ of $h_{\varphi (-,n)}$ equals ${\bf L}^{ - (n+1)d}$ times
the $n$-th coefficient of $\widehat{P}(T)$. To prove that $\widehat{P}(T)$ is a rational
power series we have to analyze how the arithmetic motivic volume of $h_{\varphi (-,n)}$
depends on $n$. To study this, we use Theorem 4.1 (quantifier elimination of Pas) to
replace the formula $\varphi (x,n)$ by a DVR-formula $\psi (x,n)$ with no quantifiers running over the
valuation ring and no quantifiers over the value group. We take an embedded resolution of
singularities $\pi: Y \longrightarrow {\bf A}_k ^d$ of the union of the loci of the
polynomials over the valuation ring, that appear in $\psi (x,n)$. Thus the pull-backs to $Y$ of
these polynomials, and the jacobian determinant of $\pi$, are locally a monomial times a
unit. Thus the pull-back of the formula $\psi (x,n)$ is easy to study, at least if one is
not scared of complicated formula in residue field variables. The key idea is to
calculate the arithmetic motivic volume of $h_{\psi (-,n)}$, by expressing it as a sum of
arithmetic motivic volumes on $h_{\mathcal{L}(Y)}$, using the change of variables formula
in Theorem 6.2. These volumes can be computed explicitly, and this yields the
rationality of $\widehat{P}(T)$.

To prove that $\widehat{P}(T)$ specializes to the Serre Poincar\'e series $P_p(T)$ for
$p\gg 0$, we repeat the above argument working with ${\bf Z}_p ^d$ instead of
${\mathcal{L}({\bf A}_k ^d)}$. The $p$-adic volume of the subset of ${\bf Z}_p ^d$ defined by
the formula $\varphi (x,n)$ equals $p^{-(n+1)d}$ times the $n$-th coefficient of
$P_p(T)$. Because of Theorem 4.1.(2), we can again replace $\varphi (x,n)$ by the formula
$\psi (x,n)$ that we obtained already above. That $p$-adic volume can be calculated explicitly by
pulling it back to the $p$-adic manifold $Y({\bf Z}_p)$, and one verifies a posteriory that
it is obtained by applying the operator $N_p$ to the arithmetic motivic volume that we
calculated above. This verification uses the last assertion in Theorem 2.1.

\section{The general setting}
\label{section 8}

We denote by $\overline{{\mathcal{M}}}$ the image of
${\rm K}_0^{mot} ({\rm Var}_k )[{\bf L}^{ -  1} ]$ in
$\widehat{\rm K}_0^{mot} ({\rm Var}_k )[{\bf L}^{ -  1} ]$, and by
$\overline{{\mathcal{M}}}_{loc}$ the localization of
$\overline{{\mathcal{M}}} \otimes {\bf Q}$ obtained by
inverting the elements ${\bf L}^i-1$, for all $i\geq 1$.
One verifies that the operator $N_p$ can be applied to any element
of $\overline{{\mathcal{M}}}_{loc}$, for $p \gg 0$, yielding a rational number. The same
holds for the Hodge-Deligne polynomial
which now belongs to ${\bf Q}(u,v)$.
By the method of section 7, we proved in \cite{Denef Loeser: Definable Sets} the following

{\bf Theorem 8.1.} \it Let $X$ be an algebraic variety
over a field $k$ of characteristic zero, let $h$ be a definable
subassignement of $h_{\mathcal{L}(X)}$, and $h_n$ a definable family of
definable subassignements of $h_{\mathcal{L}(X)}$.

\noindent (1) The motivic volume $\nu (h)$ is contained in $\overline{{\mathcal{M}}}_{loc}$.

\noindent (2) The power series $
\sum\limits_n {\nu (h_n )T} ^n \in \overline{{\mathcal{M}}}_{loc}[[T]]$ is rational, with
denominator a product of factors of the form $1-{\bf L}^{-a}T^b$, with $a,b \in{\bf N}$,
$b\neq0$.

\rm

Let $X$ be a reduced separable scheme of finite type over {\bf Z}, and let
$A=(A_p)_{p\gg0}$ be a definable family of subsets of $X({\bf Z}_p)$, meaning that on each
affine open, of a suitable finite affine covering of $X$, $A_p$ can be described by a DVR-formula
over {\bf Z}. (Here $p$ runs over all large enough primes.) To $A$ we associate in a
canonical way, its motivic volume $\nu(h_A)\in \overline{{\mathcal{M}}}_{loc}$,
in the following way: Let
$h_A$ be a definable subassignement of $h_{\mathcal{L}(X\otimes{\bf Q})}$, given by
DVR-formulas that define $A$. Because these formulas are not canonical, the subassignement
$h_A$ is not canonical. But by the Ax-Kochen-Ersov Principle (see 4.2), the set
$h_A(K)$ is canonical for each pseudo-finite field $K$ containing {\bf Q}. Hence
$\nu(h_A)\in \overline{{\mathcal{M}}}_{loc}$ is canonical,
by Theorem 6.2.(1). By the method of section 7, we proved in \cite{Denef Loeser: Definable Sets}
the following comparison result:

{\bf Theorem 8.2.} \it With the above notation, for all large enough primes $p$,
$N_p(\nu(h_A))$ equals the measure of $A_p$ with respect to the canonical measure on
$X({\bf Z}_p)$. \rm

When $X\otimes{\bf Q}$ is nonsingular and of dimension $d$, the canonical measure
on $X({\bf Z}_p)$ is defined by requiring that each fiber of the
map $X({\bf Z}_p)\rightarrow X({\bf Z}_p/p^m)$ has measure $p^{-md}$ whenever $m\gg0$.
For the definition of the canonical measure in the general case, we refer to \cite{Oesterle}.

The above theorem easily generalizes to integrals instead of measures, but this yields
little more because quite general $p$-adic integrals (such as the orbital integrals
appearing in the Langlands program) can be written as measures of the definable sets we
consider. For example the $p$-adic integral $\int {|f(x)| dx}$ on ${\bf Z}_p^d$ equals the $p$-adic
measure of
$\{ (x,t) \in {\bf Z}_p^{d + 1} :{\rm  ord}_p (f(x)) \le {\rm ord}_p (t)\}$.

\end{document}